\documentstyle[11pt,newlfont]{amsart}

\newtheorem{thm}{Theorem}[subsection]

\newtheorem{prop}[thm]{Proposition}
\newtheorem{lem}[thm]{Lemma}

\theoremstyle{definition}

\theoremstyle{remark}

\newtheorem{exa}[thm]{Example}

\newtheorem{prob}[thm]{Problem}

\numberwithin{equation}{subsection}
\newcommand{\iso}{\stackrel{\simeq}{\rightarrow}}

\newcommand{\ar}{\rightarrow}
\newcommand{\opn}{\operatorname}
\renewcommand{\mathrm}[1]{\text{\rom{\rm #1}}}

\newcommand{\nmbf}[1]{\text{\rom{\bf #1}}}
\newcommand{\bdot}{{\textstyle \cdot}}
\newcommand{\ul}{\underline}

\title{Some Remarks on Beilinson Adeles}
\author[Amnon Yekutieli]{Amnon Yekutieli*}
\address{Department of Theoretical Mathematics,
The Weizmann Institute of Science,
Rehovot 76100, ISRAEL}
\date{14 August 1995}
\email{amnon@@wisdom.weizmann.ac.il}
\thanks{* Supported by an Allon Fellowship, and incumbent of
the Anna and Maurice Boukstein Career Development Chair}

\newcommand{\lrar}[1]{\begin{picture}(50,10)(-25,-5)
\put(-25,0){\vector(1,0){50}}
\put(0,5){\makebox(0,0)[b]{\mbox{$#1$}}}
\end{picture}}

\newcommand{\ldar}[1]{\begin{picture}(10,50)(-5,-25)
\put(0,25){\vector(0,-1){50}}
\put(5,0){\mbox{$#1$}}
\end{picture}}

\newcommand{\luar}[1]{\begin{picture}(10,50)(-5,-25)
\put(0,-25){\vector(0,1){50}}
\put(5,0){\mbox{$#1$}}
\end{picture}}

\begin{document}
\maketitle

\subsection{Introduction}
In this note we consider two aspects of Beilinson adeles on schemes.

Let $X$ be a scheme of finite type over a field $k$. Given a quasi-coherent
sheaf $\cal{M}$ let
$\ul{\Bbb{A}}^{q}_{\mathrm{red}}(\cal{M})$
be the sheaf of reduced Beilinson adeles of degree $q$ (see \cite{Be},
\cite{Hr}, \cite{HY1}). It is known that
$\ul{\Bbb{A}}^{q}_{\mathrm{red}}(\cal{M}) \cong
\ul{\Bbb{A}}^{q}_{\mathrm{red}}(\cal{O}_{X})
\otimes_{\cal{O}_{X}} \cal{M}$.
For any open set $U \subset X$
\begin{equation} \label{eqn0.1}
\Gamma(U, \ul{\Bbb{A}}^{q}_{\mathrm{red}}(\cal{M})) \subset
\prod_{\xi \in S(U)^{\mathrm{red}}_{q}} \cal{M}_{\xi}
\end{equation}
where $S(U)^{\mathrm{red}}_{q}$ is the set of reduced chains of points
in $U$ of length $q$, and $\cal{M}_{\xi}$ is the Beilinson completion
of $\cal{M}$ along the chain $\xi$ (cf.\ \cite{Ye1}). For $q=0$ and $\cal{M}$
coherent one has
$\cal{M}_{(x)} = \widehat{\cal{M}}_{x}$, the $\frak{m}_{x}$-adic completion,
and (\ref{eqn0.1}) is an equality.

Let $\Omega^{\bdot}_{X/k}$ be the De Rham complex on $X$, relative to $k$.
As shown in \cite{HY1}, setting
$\cal{A}_{X}^{p, q} := \ul{\Bbb{A}}^{q}_{\mathrm{red}}(\Omega^{p}_{X/k})$
and
$\cal{A}_{X}^{i} := \bigoplus_{p+q=i} \cal{A}_{X}^{p, q}$
we get a differential graded algebra (DGA), which is quasi-isomorphic
to $\Omega^{\bdot}_{X/k}$ and is flasque.
Thus
$\mathrm{H}^{\bdot}(X, \Omega^{\bdot}_{X/k}) =
\mathrm{H}^{\bdot} \Gamma(X, \cal{A}_{X}^{\bdot})$. In particular if $X$
is smooth, we get the De Rham cohomology
$\mathrm{H}^{\bdot}_{\mathrm{DR}}(X/k)$.

More generally, let $\frak{X}$ be a formal scheme, of formally finite
type (f.f.t.) over $k$ (see \cite{Ye2}). Then applying the adelic
construction to the complete De Rham complex
$\widehat{\Omega}^{\bdot}_{\frak{X}/k}$ we get a DGA
$\cal{A}_{\frak{X}}^{\bdot}$. If $X \subset \frak{X}$ is a smooth formal
embedding (op.\ cit.) and $\opn{char} k = 0$ then
$\mathrm{H}^{\bdot} \Gamma(X, \cal{A}_{\frak{X}}^{\bdot})=
\mathrm{H}^{\bdot}_{\mathrm{DR}}(X/k)$.

There is an analogy between the sheaf $\cal{A}_{X}^{p,q}$
on a smooth $n$-dimensional variety $X$ and the sheaf of smooth
$(p,q)$-forms on a complex manifold. The coboundary operator $\mathrm{D}$ of
$\cal{A}_{X}^{\bdot}$ is defined as a sum
$\mathrm{D} := \mathrm{D}' + \mathrm{D}''$,
and
$\mathrm{D}'' : \cal{A}_{X}^{p, q} \ar \cal{A}_{X}^{p, q+1}$
plays the role of the anti-holomorphic derivative. The
map
$\int_{X} = \sum_{\xi} \opn{Res}_{\xi} : \Gamma(X, \cal{A}_{X}^{2n})
\ar k$
is the counterpart of the integral ($\opn{Res}_{\xi}$ is
the Parshin-Lomadze residue along the maximal chain $\xi$ in $X$, see
\cite{Ye1}).
This analogy to the complex manifold picture is quite
solid; for example, in \cite{HY2} there is an algebraic
proof of the Bott residue formula, which in some parts is just a
translation of the original proof of Bott to the setting of adeles.

The main purpose of this note is to examine the potential applicability of
adeles for the study of algebraic De Rham cohomology.
In \S 1 the construction of Deligne-Illusie \cite{DI} is
rewritten in terms of adeles.
In \S 2 we consider a possibility to relate adeles to Hodge theory, and
show by example its failure.

\subsection{Lifting Modulo $p^{2}$}
We interpret, in terms of adeles, the result of Deligne-Illusie
on the decomposition of the De Rham complex in characteristic $p$.
In this part we shall follow closely the ideas and notation of \cite{DI}.

Let $k$ be a perfect field of characteristic $p$.
Write $\tilde{k} := W_{2}(k)$.
Let $F_{k} : \opn{Spec} k \ar \opn{Spec} k$ be the Frobenius morphism,
i.e.\ $F^{*}_{k}(a) = a^{p}$ for $a \in k$. By pullback along $F_{k}$ we
get a scheme
$X' := X \times_{k} k$ and a finite, bijective  $k$-morphism
$F = F_{X/k} : X \ar X'$.

Assume we are given some lifting $\tilde{X}$ of $X$ to $\tilde{k}$. By this
we mean a smooth scheme $\tilde{X}$ over $\tilde{k}$ s.t.\
$X \cong \tilde{X} \times_{\tilde{k}} k$. Using the Frobenius
$F_{\tilde{k}}$ we also define a scheme
$\tilde{X}'$, and a $\tilde{k}$-morphism
$F_{\tilde{X}} : \tilde{X} \ar \tilde{X}'$.
For any point $x \in X$ the relative Frobenius homomorphism
$F^{*}_{x} : \cal{O}_{X', F(x)} \ar \cal{O}_{X, x}$
can be lifted to a $\tilde{k}$-algebra homomorphism
$\tilde{F}^{*}_{x} : \cal{O}_{\tilde{X}', F(x)} \ar
\cal{O}_{\tilde{X}, x}$
(cf.\ \cite{DI}). In view of (\ref{eqn0.1}), the collection
$\{ \tilde{F}^{*}_{x} \}_{x \in X}$
induces a homomorphism of sheaves of DG $\tilde{k}$-algebras
\[ \tilde{F}^{*} :
\ul{\Bbb{A}}^{0}_{\mathrm{red}}(\Omega^{\bdot}_{\tilde{X}' / \tilde{k}})
\ar F_{*}
\ul{\Bbb{A}}^{0}_{\mathrm{red}}(\Omega^{\bdot}_{\tilde{X} / \tilde{k}}) . \]

\begin{lem} \label{lem1}
The liftings $\{ \tilde{F}^{*}_{x} \}_{x \in X}$
determine $\cal{O}_{X'}$-linear homomorphisms
\[ f : \Omega^{1}_{X' / k} \ar F_{*} \cal{A}^{1, 0}_{X} \]
\[ h : \Omega^{1}_{X' / k} \ar F_{*} \cal{A}^{0, 1}_{X} \]
such that
\[ \mathrm{D} (f + h) = 0 . \]
\end{lem}

\begin{pf}
Let
$\nmbf{p} : \Omega^{\bdot}_{X/k} \iso
p \Omega^{\bdot}_{\tilde{X}/\tilde{k}}$
be multiplication by $p$. This extends to an
$\ul{\Bbb{A}}^{0}_{\mathrm{red}}(\cal{O}_{X})$-linear isomorphism
\[ \nmbf{p} : \cal{A}^{\bdot, 0}_{X} =
\ul{\Bbb{A}}^{0}_{\mathrm{red}}(\Omega^{\bdot}_{X/k}) \iso
p \ul{\Bbb{A}}^{0}_{\mathrm{red}}(\Omega^{\bdot}_{\tilde{X}/\tilde{k}}) . \]
Just as in \cite{DI} we get a homomorphism $f$ making the diagram
\[ \setlength{\unitlength}{0.20mm}
\begin{array}{ccc}
\Omega^{1}_{\tilde{X}'/\tilde{k}} & \lrar{\tilde{F}^{*}} &
p F_{*}
\ul{\Bbb{A}}^{0}_{\mathrm{red}}(\Omega^{1}_{\tilde{X}/\tilde{k}}) \\
\ldar{} & & \luar{\nmbf{p}} \\
\Omega^{1}_{X'/k} & \lrar{f} &
F_{*} \ul{\Bbb{A}}^{0}_{\mathrm{red}}(\Omega^{1}_{X/k})
\end{array} \]
commutative.

Next, for any chain of points $(x_{0}, x_{1})$ in $X$ and a local section
$a \in \cal{O}_{\tilde{X}'}$ we have
\[ \mathrm{D}'' \tilde{F}^{*}(a) =
\tilde{F}^{*}_{x_{0}}(a) - \tilde{F}^{*}_{x_{1}}(a) \in
p \cal{O}_{\tilde{X}, (x_{0}, x_{1})} . \]
Therefore
\[ \mathrm{D}'' \tilde{F}^{*} : \cal{O}_{\tilde{X}'} \ar
p F_{*} \ul{\Bbb{A}}^{1}_{\mathrm{red}}(\cal{O}_{\tilde{X}}) \]
is a derivation which kills $p \cal{O}_{\tilde{X}'}$, and we get an
$\cal{O}_{X'}$-linear homomorphism $h$ s.t.\ the diagram
\[ \setlength{\unitlength}{0.20mm}
\begin{array}{ccc}
\cal{O}_{\tilde{X}'} & \lrar{\mathrm{D}'' \tilde{F}^{*}} &
F_{*} p
\ul{\Bbb{A}}^{1}_{\mathrm{red}}(\cal{O}_{\tilde{X}}) \\
\ldar{\mathrm{d}} & & \luar{\nmbf{p}} \\
\Omega^{1}_{X'/k} & \lrar{h} &
F_{*} \ul{\Bbb{A}}^{1}_{\mathrm{red}}(\cal{O}_{X})
\end{array} \]
commutes.

Reinterpreting the calculations of \cite{DI} in terms of adeles we see that
the following hold: for each point $x_{0} \in X$,
$\mathrm{D}' f = 0$ in $\Omega^{1}_{X/k, (x_{0})}$;
for each chain $(x_{0}, x_{1}, x_{2})$ in $X$,
$\mathrm{D}'' h = 0$ in $\cal{O}_{X, (x_{0}, x_{1}, x_{2})}$;
lastly, for each chain $(x_{0}, x_{1})$,
$\mathrm{D}'' f = - \mathrm{D}' h$ in $\Omega^{1}_{X/k, (x_{0}, x_{1})}$.
This implies that on the level of sheaves $\mathrm{D}(f + h) = 0$.
\end{pf}

\begin{prop} \label{prop1}
The liftings $\{ \tilde{F}^{*}_{x} \}_{x \in X}$
determine an $\cal{O}_{X'}$-linear homomorphism of complexes
\[ \psi_{\tilde{X}} : \bigoplus_{i = 0}^{n}
\ul{\Bbb{A}}^{\bdot}_{\mathrm{red}}(\Omega^{i}_{X' / k})[-i] \ar
F_{*} \cal{A}^{\bdot}_{X} \]
making the diagram
\begin{equation} \label{eqn6}
\setlength{\unitlength}{0.25mm}
\begin{array}{ccc}
\bigoplus_{i} \Omega^{i}_{X' / k} & \lrar{C^{-1}} &
F_{*} \mathrm{H}^{\bdot} \Omega^{\bdot}_{X / k} \\
\ldar{} & & \ldar{} \\
\bigoplus_{i} \mathrm{H}^{\bdot}
\ul{\Bbb{A}}^{\bdot}_{\mathrm{red}}(\Omega^{i}_{X' / k})[-i] &
\lrar{\mathrm{H}^{\bdot}(\psi_{\tilde{X}})} &
F_{*} \mathrm{H}^{\bdot}  \cal{A}^{\bdot}_{X}
\end{array}
\end{equation}
commute. Here $C^{-1}$ is the Cartier operation, and the vertical arrows
are the canonical isomorphisms. Therefore $\psi_{\tilde{X}}$ is a
quasi-isomorphism.
\end{prop}

\begin{pf}
Since
\[ \ul{\Bbb{A}}^{j}_{\mathrm{red}}(\Omega^{i}_{X' / k}) \cong
\ul{\Bbb{A}}^{j}_{\mathrm{red}}(\cal{O}_{X'}) \otimes_{\cal{O}_{X'}}
\Omega^{i}_{X' / k} , \]
and since
$F^{*} : \ul{\Bbb{A}}^{\bdot}_{\mathrm{red}}(\cal{O}_{X'}) \ar
F_{*} \ul{\Bbb{A}}^{\bdot}_{\mathrm{red}}(\cal{O}_{X})$
commutes with $\mathrm{D}''$ and is killed by $\mathrm{D}'$
(i.e.\ $\mathrm{D}' F^{*} = 0$)
it suffices to define $\cal{O}_{X'}$-linear homomorphisms
$\psi_{\tilde{X}}^{i} : \Omega^{i}_{X' / k} \ar F_{*} \cal{A}^{i}_{X}$
s.t.\ $\mathrm{D} \psi_{\tilde{X}}^{i} = 0$.
Define
$\psi_{\tilde{X}}^{0} := F^{*}$, and
$\psi_{\tilde{X}}^{1} := f+h$ as in Lemma \ref{lem1}. For
$1 \leq i \leq n$ let
$\nmbf{a} : \Omega^{i}_{X' / k} \ar (\Omega^{1}_{X' / k})^{\otimes i}$
be the anti-symmetrizing operator (this makes sense since $n < p$;
cf.\ \cite{DI}), and define $\psi_{\tilde{X}}^{i}$ by
\[ \setlength{\unitlength}{0.30mm}
\begin{array}{ccc}
(\Omega^{1}_{X' / k})^{\otimes i} &
\lrar{(\psi_{\tilde{X}}^{1})^{\otimes i}} &
(F_{*} \cal{A}^{1}_{X})^{\otimes i} \\
\luar{\nmbf{a}} & & \ldar{\mathrm{product}} \\
\Omega^{i}_{X' / k} & \lrar{\psi_{\tilde{X}}^{i}} &
F_{*} \cal{A}^{i}_{X}
\end{array} \]

Let $a \in \cal{O}_{\tilde{X}}$
be a local section, with corresponding pullback
$a \otimes 1 \in \cal{O}_{\tilde{X}'}$, and with image
$a_{0} \in \cal{O}_{X}$. Then according to the calculations in \cite{DI},
we have
$\tilde{F}^{*}(a \otimes 1) = a^{p} + \nmbf{p} u$
for some local section
$u \in \ul{\Bbb{A}}^{0}_{\mathrm{red}}(\cal{O}_{X})$. Therefore
$f (\mathrm{d} a_{0} \otimes 1) = a_{0}^{p-1} \mathrm{d} a_{0} +
\mathrm{D}' u$ and
$h (\mathrm{d} a_{0} \otimes 1) = \mathrm{D}'' u$,
so
\[ \psi_{\tilde{X}} (\mathrm{d} a_{0} \otimes 1) =
a_{0}^{p-1} \mathrm{d} a_{0} + \mathrm{D} u . \]
This means that
\[ \mathrm{H}^{1} (\psi_{\tilde{X}}) = C^{-1} : \Omega^{1}_{X' / k} \iso
F_{*} \mathrm{H}^{1} \Omega^{\bdot}_{X / k} \cong
F_{*} \mathrm{H}^{1} \cal{A}^{\bdot}_{X} . \]
Clearly in degree $0$,
$\mathrm{H}^{0} (\psi_{\tilde{X}}) = F^{*} = C^{-1}$.
Since the vertical arrows in diagram (\ref{eqn6}) are isomorphisms of
(sheaves of) graded algebras, it follows that
$\mathrm{H}^{\bdot} \cal{A}^{\bdot}_{X}$ is a graded-commutative algebra,
and therefore
\[ \mathrm{H}^{\bdot}(\psi_{\tilde{X}}) :
\bigoplus_{i} \Omega^{i}_{X' / k} \ar F_{*} \mathrm{H}^{\bdot}
\cal{A}^{\bdot}_{X} \]
is a homomorphism of graded algebras. But then
$\mathrm{H}^{\bdot}(\psi_{\tilde{X}}) = C^{-1}$ in all degrees, and it's
an isomorphism.
\end{pf}

Of course in the derived category the map
$\psi_{\tilde{X}}$ is independent of the choices of Frobenius liftings.

\subsection{A Hodge-type Decomposition?}
The second aspect is a naive attempt to use adeles for a Hodge-type
decomposition of De Rham cohomology. Suppose  $\opn{char} k = 0$  and
$X$ is smooth over $k$, of dimension $n$. For any $0 \leq p, q \leq n$
define a canonical subspace
\begin{equation}
\mathrm{H}^{p,q} :=
\frac{ \Gamma(X, \cal{A}^{p,q}_{X}) \cap \opn{Ker} D }
{ \Gamma(X, \cal{A}^{p,q}_{X}) \cap \opn{Im} D }
\subset \mathrm{H}^{p+q}_{\mathrm{DR}}(X/k)
\end{equation}
(cf.\ \cite{GH} p.\ 116).
Since the sheaves $\cal{A}^{p, q}_{X}$ imitate the Dolbeault sheaves on
a complex manifold so nicely, one can imagine that
\[ \mathrm{H}^{i}_{\mathrm{DR}}(X/k) =
\bigoplus_{p+q=i} \mathrm{H}^{p, q} \]
if $X$ is proper.
Yet this is {\em false}, as can be seen from the example below.

What we get is a serious breakdown in the analogy to smooth forms on
a complex manifold. I should mention that even in \cite{HY2} there
was a  breakdown in this analogy; there it was not possible to define
a connection on the adelic sections of a vector bundle, and hence
an auxiliary algebraic device, the sheaf
$\tilde{\cal{A}}_{X}^{\bdot}$ of Thom-Sullivan adeles,
had to be introduced.

\begin{prob}
Is it true that for $X$ smooth, the filtration on
$\cal{A}_{X}^{\bdot}$ by the subcomplexes
$\cal{A}_{X}^{\bdot, \geq q}$
induces the coniveau filtration on
$\mathrm{H}^{\bdot}_{\mathrm{DR}}(X/k)$?
\end{prob}

\begin{exa} \label{exa1}
Suppose $k$ is algebraically closed and $X$ is an elliptic
curve. Then
$\opn{dim} \mathrm{H}^{1}_{\mathrm{DR}}(X/k) = 2$.
Consider the nondegenerate pairing on $\mathrm{H}^{1}_{\mathrm{DR}}(X/k)$
given by
\[ \langle [\alpha] , [\beta] \rangle =
\int_{X} [\alpha] \smile [\beta] =
\sum_{\xi} \opn{Res}_{\xi} (\alpha \cdot \beta) \]
for adeles $\alpha, \beta \in \cal{A}^{1}_{X}$. We see that
$\langle \mathrm{H}^{1, 0} , \mathrm{H}^{1, 0} \rangle =
\langle \mathrm{H}^{0, 1} , \mathrm{H}^{0, 1} \rangle  = 0$.
Therefore if
$\mathrm{H}^{1}_{\mathrm{DR}}(X) = \mathrm{H}^{1, 0} + \mathrm{H}^{0, 1}$,
then
$\opn{dim} \mathrm{H}^{1, 0} = \opn{dim} \mathrm{H}^{0, 1} = 1$.
It is easy to find
$0 \neq [\alpha] \in \mathrm{H}^{1, 0}$;
take any $0 \neq [\alpha] \in \Gamma(X, \Omega^{1}_{X/k})$.
On the other hand an adele
\[ \beta = (b_{(\mathrm{gen}, x)}) \in \Gamma(X, \cal{A}^{0, 1}_{X}) =
\Bbb{A}^{1}_{\mathrm{red}}(X, \cal{O}_{X}) \]
(where $x$ runs over the set $X_{0}$ of closed points, and $\mathrm{gen}$
is the generic point) satisfies $\mathrm{D} \beta = 0$ iff
$\mathrm{d} b_{(\mathrm{gen}, x)} = 0$ for every $x$. This forces
$b_{(\mathrm{gen}, x)} \in k$. But taking
$b = (b_{(\mathrm{gen})}, b_{(x)}) \in
\Bbb{A}^{0}_{\mathrm{red}}(X, \cal{O}_{X})$, with
$b_{(\mathrm{gen})} = 0$, $b_{(x)} = b_{(\mathrm{gen}, x)}$
we get $\beta = \mathrm{D} b$. Hence
$\mathrm{H}^{0, 1} = 0$.
\end{exa}

\begin{prob} \label{prob2}
For $\alpha$ as above find explicitly a cocycle
$\beta \in \Gamma(X, \cal{A}^{1}_{X})$ s.t.\
$\langle \alpha , \beta \rangle = 1$.
\end{prob}

The best I can do is:

\begin{prop} \label{prop2}
Suppose $X$ is a smooth proper curve and $k$ is algebraic\-ally closed.
Let $\alpha_{(\mathrm{gen})} \in \Omega^{1}_{k(X) / k}$ be a differential
of the $2$-nd kind, namely
$\opn{Res}_{(\mathrm{gen}, x)} \alpha_{(\mathrm{gen})} = 0$
for every $x \in X_{0}$.
Then it defines a cocycle
$\alpha \in \Gamma(X, \cal{A}^{1}_{X})$
whose component at $(\mathrm{gen})$ is $\alpha_{(\mathrm{gen})}$.
Every cohomology class in $\mathrm{H}^{1}_{\mathrm{DR}}(X / k)$
is gotten in this way. The Hodge filtration is induced by the differentials
of the $1$-st kind.
\end{prop}

\begin{pf}
The adele $\alpha$ will be given by its bihomogeneous components,
$\alpha = \alpha^{1,0} + \alpha^{0,1}$. We set
$\alpha^{1,0} := (\alpha_{(\mathrm{gen})}, \alpha_{(x)})$
where for $x \in X_{0}$, $\alpha_{(x)} = 0$.
Since
$\opn{Res}_{(\mathrm{gen}, x)} \alpha_{(\mathrm{gen})} = 0$
there is some
$a_{(\mathrm{gen}, x)} \in k(X)_{(\mathrm{gen}, x)}$
(unique up to adding a constant) s.t.\
$\mathrm{d} a_{(\mathrm{gen}, x)} = \alpha_{(\mathrm{gen})}$.
Set
$\alpha^{0,1} := ( a_{(\mathrm{gen}, x)} )$.
Then $\alpha$ is evidently a cocycle.

If $\alpha_{(\mathrm{gen})}$ is of the $1$-st kind then actually we get
$a_{(\mathrm{gen}, x)} \in \cal{O}_{X, (x)}$; call this element also
$a_{(x)}$.
So we can define an adele
$\tilde{\alpha} = \tilde{\alpha}^{1,0} + \tilde{\alpha}^{0,1}$
with
$\tilde{\alpha}^{1,0} := (\alpha_{(\mathrm{gen})},  \mathrm{d} a_{(x)} )$
and
$\tilde{\alpha}^{0,1} := 0$. We get a cocycle (cohomologous to $\alpha$),
and conversely any
cocycle in $\Gamma(X, \cal{A}^{1,0}_{X})$ looks like this.

Consider the niveau spectral sequence of De Rham homology (cf.\ \cite{Ye3}).
A comparison of dimensions shows that this degenerates at the $E_{2}$
term. Also the niveau filtration on $\mathrm{H}_{1}^{\mathrm{DR}}(X / k)$
is trivial. Hence we get
\[ \begin{array}{rcl}
\mathrm{H}_{1}^{\mathrm{DR}}(X / k) & = &
\opn{Ker} \left( \mathrm{H}^{1} \Omega^{\bdot}_{k(X) / k} \ar
\bigoplus_{x \in X_{0}} k \right) \\
& \cong & (\text{forms of the $2$-nd kind}) / (\text{exact forms}) .
\end{array} \]

Now the map
$\mathrm{H}^{1}_{\mathrm{DR}}(X / k) \ar
\mathrm{H}_{1}^{\mathrm{DR}}(X / k)$,
$[\alpha] \mapsto [\alpha] \frown [X] = \pm [\mathrm{C}_{X} \cdot \alpha]$
is bijective. A direct inspection reveals that the adele
$\alpha = \alpha^{1,0} + \alpha^{0,1}$
is sent to the differential of the second kind
$\alpha_{(\mathrm{gen})} \in \Omega^{1}_{k(X)/k}$.
\end{pf}


\end{document}